\documentclass[12pt]{article}
\usepackage{amssymb, amsmath, amsthm, amscd}
\usepackage[pdftex]{graphics}
\usepackage[utf8]{inputenc}
\usepackage[all,cmtip]{xy}
\usepackage{bbm}
\usepackage{enumitem}
\usepackage{setspace}
\usepackage{mathdots}

\usepackage[mathscr]{euscript}
\usepackage[colorlinks=true]{hyperref}
\hypersetup{colorlinks   = true}
\hypersetup{linkcolor=blue}
\usepackage{tikz}
\usetikzlibrary{cd}

\begin{document}

\mathchardef\mhyphen="2D
\newtheorem{The}{Theorem}[section]
\newtheorem{Lem}[The]{Lemma}
\newtheorem{Prop}[The]{Proposition}
\newtheorem{Cor}[The]{Corollary}
\newtheorem{Rem}[The]{Remark}
\newtheorem{Obs}[The]{Observation}
\newtheorem{SConj}[The]{Standard Conjecture}
\newtheorem{Titre}[The]{\!\!\!\! }
\newtheorem{Conj}[The]{Conjecture}
\newtheorem{Question}[The]{Question}
\newtheorem{Prob}[The]{Problem}
\newtheorem{Def}[The]{Definition}
\newtheorem{Not}[The]{Notation}
\newtheorem{Claim}[The]{Claim}
\newtheorem{Conc}[The]{Conclusion}
\newtheorem{Ex}[The]{Example}
\newtheorem{Fact}[The]{Fact}
\newtheorem{Formula}[The]{Formula}
\newtheorem{Formulae}[The]{Formulae}
\newtheorem{The-Def}[The]{Theorem and Definition}
\newtheorem{Prop-Def}[The]{Proposition and Definition}
\newtheorem{Lem-Def}[The]{Lemma and Definition}
\newtheorem{Cor-Def}[The]{Corollary and Definition}
\newtheorem{Conc-Def}[The]{Conclusion and Definition}
\newtheorem{Terminology}[The]{Note on terminology}
\newcommand{\C}{\mathbb{C}}
\newcommand{\R}{\mathbb{R}}
\newcommand{\N}{\mathbb{N}}
\newcommand{\Z}{\mathbb{Z}}
\newcommand{\Q}{\mathbb{Q}}
\newcommand{\Proj}{\mathbb{P}}
\newcommand{\Rc}{\mathcal{R}}
\newcommand{\Oc}{\mathcal{O}}
\newcommand{\Vc}{\mathcal{V}}
\newcommand{\Id}{\operatorname{Id}}
\newcommand{\pr}{\operatorname{pr}}
\newcommand{\rk}{\operatorname{rk}}
\newcommand{\del}{\partial}
\newcommand{\delbar}{\bar{\partial}}
\newcommand{\Cdot}{{\raisebox{-0.7ex}[0pt][0pt]{\scalebox{2.0}{$\cdot$}}}}
\newcommand\nilm{\Gamma\backslash G}
\newcommand\frg{{\mathfrak g}}
\newcommand{\fg}{\mathfrak g}
\newcommand{\Oh}{\mathcal{O}}
\newcommand{\Kur}{\operatorname{Kur}}
\newcommand\gc{\frg_\mathbb{C}}
\newcommand\jonas[1]{{\textcolor{green}{#1}}}
\newcommand\luis[1]{{\textcolor{red}{#1}}}
\newcommand\dan[1]{{\textcolor{blue}{#1}}}

\begin{center}

{\Large\bf A Comparative Study between K\"ahler and Non-K\"ahler Hyperbolicity}

\end{center}

\begin{center}
{\large Samir Marouani 
}

\end{center}

\vspace{1ex}

\noindent{\small{\bf Abstract.} In this note, we establish a connection between SKT  and balanced hyperbolicity, highlighting their relationship with Kähler hyperbolicity in the sense of Gromov.

\vspace{2ex}
\section{Introduction}\label{section:Introduction}
In this paper, we continue the study of compact complex {\it SKT hyperbolic} manifolds that we introduced very recently in [Mar23] as generalisations in the possibly non-K\"ahler context of the classical notion of {\it K\"ahler hyperbolic} (in the sense of Gromov) manifolds. 

\vspace{2ex}
According to [Gro91], a $C^\infty$ $k$-form $\alpha$ on $X$ is said to be $\widetilde{d}(\mbox{bounded})$ with respect to $\omega$ if $\pi_X^\star\alpha = d\beta$ on $\widetilde{X}$ for some $C^\infty$ $(k-1)$-form $\beta$ on $\widetilde{X}$ that is bounded w.r.t. $\widetilde\omega$.
Now, recall that a compact complex manifold $X$ is said to be {\it K\"ahler hyperbolic} in the sense of Gromov (see [Gro91]) if there exists a K\"ahler metric $\omega$ on $X$ (i.e. a Hermitian metric $\omega$ with $d\omega=0$) such that $\omega$ is $\widetilde{d}(\mbox{bounded})$ with respect to itself. In [MP22, Definition 2.1], we introduced the following 1-codimensional analogue of this:\\

\noindent {\it An $n$-dimensional compact complex manifold $X$ is said to be {\bf balanced hyperbolic} if there exists a balanced metric $\omega$ on $X$  (i.e. a Hermitian metric $\omega$ with $d\omega^{n-1}=0$) such that  $\omega^{n-1}$ is $\widetilde{d}(\mbox{bounded})$ with respect to $\omega$.} 
\noindent Any such metric $\omega$ is called a {\it balanced hyperbolic} metric. Meanwhile, it is well known that balanced metrics  may exist on certain compact complex $n$-dimensional manifolds $X$ such that $\omega^{n-1}$ is even $d$-exact on $X$. These manifolds are called {\it degenerate balanced}. There is no analogue of this phenomenon in the K\"ahler setting.

\vspace{2ex}
In [Mar23] we generalize the notions of $\widetilde{d}(\mbox{bounded})$ form  and the K\"ahler hyperbolicity  by introducing the following:

\vspace{2ex}
(i)\ \noindent {\it
A $C^\infty$ $k$-form $\phi$ on a complex manifold $X$ is said to be $\widetilde{(\partial+\bar\partial)}$-bounded with respect to $\omega$ if $\pi_X^\star\phi = \partial\alpha+\bar\partial\beta$ on $\widetilde{X}$ for some $C^\infty$ $(k-1)$-forms $\alpha$ and $\beta$ on $\widetilde{X}$ that are bounded w.r.t. $\widetilde\omega$.}
    
\vspace{1ex}
(ii)\ \noindent {\it  A Hermitian metric $\omega$ on a compact complex manifold $X$  is said to be {\bf SKT hyperbolic} if $\omega$ is SKT and  $(\widetilde{\partial+\bar\partial})-\mbox{bounded}$ with respect to $\omega$.}\\
 The manifold $X$ is said to be SKT hyperbolic if it carries a {\it SKT hyperbolic metric}.
 
 \vspace{2ex}
 
On the other hand, $X$ is said to be {\it Kobayashi hyperbolic} (see e.g. [Kob70]) if the Kobayashi pseudo-distance on $X$ is actually a distance. By Brody's Theorem 4.1. in [Bro78], this is equivalent to the non-existence of entire curves in $X$, namely the non-existence of non-constant holomorphic maps $f:\C\longrightarrow X$. This latter property has come to be known as the {\it Brody hyperbolicity} of $X$. Thus, a compact manifold $X$ is Brody hyperbolic if and only if it is Kobayashi hyperbolic. (The equivalence is known to fail when $X$ is non-compact.) 

   \vspace{2ex}
   
The implications among these notions are summed up in the following diagram. (See [Mar23, Theorem 1.1] .)

\vspace{3ex}
$\begin{array}{lll} X \hspace{1ex} \mbox{is {\bf K\"ahler hyperbolic}} \implies  X \hspace{1ex} \mbox{is {\bf SKT hyperbolic}} \implies  X \hspace{1ex} \mbox{is {\bf Kobayashi hyperbolic}} \\
 \rotatebox{-90}{$\implies$} \\
 X \hspace{1ex} \mbox{is {\bf balanced hyperbolic}}  \\
 \rotatebox{90}{$\implies$} &  &   \\
 X \hspace{1ex} \mbox{is {\bf degenerate balanced}}. & &\end{array}$ 

\vspace{2ex}
On compact complex surfaces, an application of the Stokes Theorem shows that any SKT metric gives rise to a non-zero class in Aeppli cohomology, i.e., ($0\neq [\omega]_A\in H^{1,1}_{A}(X,\C)$ (see, e.g., [HL83, Prop. 37]). Very recently Yau, Zhao and Zheng prove that if $(M,J,g,\omega)$ is a compact SKL (Strominger K\"ahler-like) non-K\"ahler manifold, then $\omega$ gives rise to a non-zero class in Aeppli cohomology, i.e.,
$0\neq [\omega]_A\in H^{1,1}_{A}(X,\C)$ (see [YZZ22, Thm. 1]). Our first main result is the following:
 \begin{The}
      Let $(X,\omega)$ be an SKT compact complex $n$-dimensional manifold. Then, any {\bf SKT} metric $\omega$ gives rise to a non-zero class in Aeppli cohomology $(i.e. [\omega]_A\neq 0 )$.
 \end{The}
Therefore, we deduce the existence of hyperbolic manifolds in the sense of Kobayashi that do not exhibit SKT hyperbolicity. Now recall the following conjecture 
\begin{Conj}(Kobayashi, see e.g. [CY18, Conjecture 2.8] or Lang's survey cited therein)\label{Conj:Kobayashi}

  If $X$ is a {\bf Kobayashi hyperbolic} compact complex manifold, its canonical bundle $K_X$ is {\bf ample}.

\end{Conj} 
Therefore, as a perspective we ask the following question:
\begin{Question}
Is it true that if $X$ is a compact complex manifold that is SKT hyperbolic, then its canonical bundle $K_X$ is ample?
 \end{Question}
  \vspace{1ex}
 Our second main result is to give the link between K\"ahler and non-K\"ahler hyperbolicity; Let recall the following

\begin{Conj}(Fino-Vezzoni)\label{Conj:Fino-Vezzoni}

  Every compact complex manifold admitting both an SKT metric and a balanced metri is necessarily K\"ahler

\end{Conj}
Therefore, our approach is as follows
\begin{The}
    Let $(X,\omega)$, be a compact complex manifold with $\mbox{dim}_\C X=n$ such that $\omega$ is both SKT hyperbolic and balanced hyperbolic. Let $\pi:\widetilde{X}\longrightarrow X$ be the universal cover of $X$ and $\widetilde\omega:=\pi^\star\omega$ the lift to $\widetilde{X}$ of  $\omega$ on $X$ and $\Gamma$ be a  smooth bounded $(2n-3)$-form such that $\widetilde{\omega}^{n-1}=d\Gamma$. Then, the following equivalence holds:\\
    
     (i)\ 
       $ \quad \partial\Gamma^{n-1,n-2}=-\bar\partial\Gamma^{n,n-3}=\overline{\bar\partial\Gamma^{n-2,n-1}}=-\overline{\partial\Gamma^{n-3,n}}\in\mbox{L}^2(\widetilde{X}) $
       
        (ii)\ $\quad\omega$  is a K\"ahler hyperbolic metric in the sense of Gromov .\\

     Where $\Gamma^{n-k,n-3+k}$  is the  $(n-k, n-3+k)$-component of $\Gamma$, $k\in\{0,1,2,3\}$.
\end{The}
\section{SKT, balanced and K\"ahler hyperbolicity}

We begin our paper by presenting the proof of the first main result that follows.
 \begin{The}
      Let $(X,\omega)$ be an SKT compact complex $n$-dimensional manifold. Then, any {\bf SKT} metric $\omega$ gives rise to a non-zero class in Aeppli cohomology $(i.e. [\omega]_A\neq 0 )$.
 \end{The}
 \noindent {\it Proof.}
  Suppose there is a $(0,1)$-form $\alpha$ and a $(1,0)$-form $\beta$ such that $\omega=\partial\alpha+\bar\partial\beta$. 
We get \begin{eqnarray*} \big(d(\alpha+\beta)\big)^n =(\omega+\bar\partial\alpha+\partial\beta)^n&=&\sum_{k=0}^n\sum_{j=0}^k C_n^kC_j^k(\bar\partial\alpha)^j\wedge(\partial\beta)^{k-j}\wedge\omega^{n-k}\\ &=& \sum_{\{k=2j\}}\mathcal{C}_j(\bar\partial\alpha)^j\wedge(\partial\beta)^{j}\wedge\omega^{n-2j}. \end{eqnarray*} To show that $\sum_{\{k=2j\}}\mathcal{C}_j(\bar\partial\alpha)^j\wedge(\partial\beta)^{j}\wedge\omega^{n-2j}>0$, it suffices to check that the real form $(\bar\partial\alpha)^j\wedge(\partial\beta)^{j}\wedge\omega^{n-2j}$ is
weakly (semi)-positive at every point of $X$. (Recall that $\bar\partial\alpha$ is the conjugate of $\partial\beta$.) To this end,
note that the $(2j, 2j)$-form $(\bar\partial\alpha)^j\wedge(\partial\beta)^{j}$
is weakly semi-positive as the wedge product of a $(2j, 0)$-
form and its conjugate (see [Dem97, Chapter III, Example 1.2]). Therefore, the $(n, n)$-form is (semi)-positive since the product of a weakly (semi)-positive form and a strongly (semi)-positive form is weakly (semi)-positive and $\omega$ is strongly positive (see [Dem97,
Chapter III, Proposition 1.11]). (Recall that in bidegrees $(0, 1)$, $(1,1)$ , $(n-1,n-1)$ and $(n, n)$, the notions of
weak and strong positivity coincide).
By Stokes, we get $$0=\int_X d\big((\alpha+\beta)\wedge(d(\alpha+\beta)^{n-1}\big)=\int_X\sum_{\{k=2j\}}\mathcal{C}_j(\bar\partial\alpha)^j\wedge(\partial\beta)^{j}\wedge\omega^{n-2j}>0.$$ Consequently \begin{eqnarray}\label{non_0_A_classe}  [\omega]_A\neq 0 .\end{eqnarray}
\hfill $\Box$

We establish the following corollary:
\begin{Cor}
    The implication from SKT hyperbolicity to Kobayashi hyperbolicity is strictly one-sided.
\end{Cor}
\noindent {\it Proof.}
A very generic surfece $X\in \Proj^3$ of high degree (see [DE00] for degree $\geq 21$) is Kobayshi hyperbolic, but can not SKT hyperbolic since it is simply connected.
\hfill $\Box$

  \vspace{3ex}
  
Throughout the remainder of the text, we adopt the notation: $$\gamma_p:=\frac{\gamma^p}{p!}.$$
On a complex manifold $X$ with $\mbox{dim}_\C X=n$, we will often use the following standard formula (cf. e.g. [Voi02, Proposition 6.29, p. 150]) for the Hodge star operator $\star = \star_\omega$ of any Hermitian metric $\omega$ applied to $\omega$-{\it primitive} forms $v$ of arbitrary bidegree $(p, \, q)$: \begin{eqnarray}\label{eqn:prim-form-star-formula-gen}\star\, v = (-1)^{k(k+1)/2}\, i^{p-q}\,\omega_{n-p-q}\wedge v, \hspace{2ex} \mbox{where}\,\, k:=p+q.\end{eqnarray} Recall that, for any $k=0,1, \dots , n$, a $k$-form $v$ is said to be $(\omega)$-{\it primitive} if $\omega_{n-k+1}\wedge v=0$ and that this condition is equivalent to $\Lambda_\omega v=0$, where $\Lambda_\omega$ is the adjoint of the operator $\omega\wedge\cdot$ (of multiplication by $\omega$) w.r.t. the pointwise inner product $\langle\,\,,\,\,\rangle_\omega$ defined by $\omega$.
\begin{Lem}\label{Lem:1-forms_Delta-harm} Let $\omega$ be a Hermitian metric on a complex manifold $X$ with $\mbox{dim}_\C X=n$. Fix a primitive form $\phi \in C^\infty_{p,q}(X,\,\C)$ with $p+q=k\leq n$.

  \vspace{3ex}

If $\omega$ is {\bf K\"ahler} and $\partial\phi=\bar\partial\phi=0$, then 
\begin{eqnarray}\label{complete_case}\Delta(\omega_{n-k}\wedge\phi)=0\end{eqnarray}
 where $\Delta = \Delta_\omega = dd^\star + d^\star d$ is the $d$-Laplacian induced by $\omega$.

\end{Lem} 
\noindent {\it Proof.} 
Let $\phi$ be a primitive form, the standard formula (\ref{eqn:prim-form-star-formula-gen}) we get: $\star\, \phi = (-1)^\frac{k(k+1)}{2}i^{p-q}\omega_{n-k}\wedge \phi$, hence $\star\,(\omega_{n-k}\wedge \phi) = -(-1)^\frac{k(k+1)}{2}i^{q-p}\phi$. Meanwhile, $d^\star = -\star d\star$, so applying $-\star d$ to the previous identity, we get \begin{eqnarray*}d^\star(\omega_{n-k}\wedge \phi)=(-1)^\frac{k(k+1)}{2}i^{q-p}(\star\partial \phi +\star\bar\partial \phi).\end{eqnarray*} 
By assumption, we have $\partial\phi=\bar\partial\phi=0$, then, since $\omega$ is K\"ahler we get: \begin{eqnarray*}\Delta(\omega_{n-k}\wedge \phi) &=& 0 .\end{eqnarray*}
\hfill $\Box$
 \begin{Prop}\label{Equi_SKT_Kahler}
      Let $\pi:\widetilde{X}\longrightarrow X$ be the universal cover of a compact complex $n$-dimensional K\"ahler manifold $X$, with $n>2$ and $\widetilde\omega:=\pi^\star\omega$ the lift to $\widetilde{X}$ of an SKT hyperbolic metric $\omega$ on $X$ such that $\widetilde{\omega}=\partial\alpha+\bar\partial\beta$ for some smooth bounded $(0,1)$-form $\alpha$ and $(1,0)$-form $\beta$. Then, $\omega$ is K\"ahler hyperbolic if and only if $\bar\partial\alpha\in$ L$^{2}(\widetilde{X})$.
 \end{Prop}
 \noindent {\it Proof.}
 We will show that $\bar\partial\alpha\in$ L$^{2}(\widetilde{X})$ is equivalent to $\bar\partial\alpha=0$ and $\partial\beta=0$ hence to $\widetilde\omega=d(\alpha+\beta)$, therefore the K\"ahler hyperbolicity of $\omega$. Since all $(0,p)$-forms are primitive  by Lemma \ref{Lem:1-forms_Delta-harm}, we get a  $\Delta$-harmonic form namely $\widetilde\omega^{n-2}\wedge\bar\partial\alpha$. Since $\bar\partial\alpha\in$ L$^{2}(\widetilde{X})$ by assumption is so $\widetilde\omega^{n-2}\wedge\bar\partial\alpha$, then using Theorem 3.20 in [Mar23] we get $\widetilde\omega^{n-2}\wedge\bar\partial\alpha=0$. And since $\widetilde\omega$ is reel, we have   $\overline{\bar\partial\alpha}=\partial\bar\alpha=\partial\beta$. On the other hand, the pointwise {\it Lefschetz map}: \begin{equation*}L_{\widetilde{\omega}_{n-2}}:\Lambda^2T^\star X\longrightarrow\Lambda^{2n-2}T^\star X, \hspace{3ex} \varphi\longmapsto\psi:=\widetilde{\omega}_{n-2}\wedge\varphi,\end{equation*} is {\it bijective}.
By replacing $\varphi$ with $\bar\partial\alpha$ on one hand and $\varphi$ with $\partial\beta $ on the other hand, we obtain our desired result.
 \hfill $\Box$
 
 \vspace{1ex}
 
 Now we proceed to prove the following main result
\begin{The}
    Let $(X,\omega)$, be a compact complex manifold with $\mbox{dim}_\C X=n$ such that $\omega$ is both SKT hyperbolic and balanced hyperbolic. Let $\pi:\widetilde{X}\longrightarrow X$ be the universal cover of $X$ and $\widetilde\omega:=\pi^\star\omega$ the lift to $\widetilde{X}$ of  $\omega$ on $X$ and $\Gamma$ be a  smooth bounded $(2n-3)$-form such that $\widetilde{\omega}^{n-1}=d\Gamma$. Then, the following equivalence holds:\\
    
     (i)\ 
       $ \quad \partial\Gamma^{n-1,n-2}=-\bar\partial\Gamma^{n,n-3}=\overline{\bar\partial\Gamma^{n-2,n-1}}=-\overline{\partial\Gamma^{n-3,n}}\in\mbox{L}^2(\widetilde{X}) $
       
        (ii)\ $\quad\omega$  is a K\"ahler hyperbolic metric in the sense of Gromov .\\

     Where $\Gamma^{n-k,n-3+k}$  is the  $(n-k, n-3+k)$-component of $\Gamma$, $k\in\{0,1,2,3\}$.
\end{The}
\noindent {\it Proof.} Let  $\pi_X:\widetilde{X}\longrightarrow X$ be the universal cover of $X$.
The SKT hyperbolicity assumption on $\omega$ translates to the following:\begin{eqnarray}\label{skt_hyper_1}
\pi_X^\star\omega = \partial\alpha+\bar\partial\beta  \hspace{3ex} \mbox{on}\hspace{1ex} \widetilde{X},\end{eqnarray} where $\alpha$ and $\beta$ are a smooth, $\widetilde\omega$-bounded $(0,1)$-form and $(1,0)$-form  on $\widetilde{X}$ and $\widetilde\omega=\pi_X^\star\omega$ is the lift of the metric $\omega$ to $\widetilde{X}$.

On the other hand, the balanced hyperbolicity assumption on $\omega$ translates to the following propertie:
\begin{eqnarray}\label{balanced_hyper_1}\pi_X^\star\omega^{n-1} = d\Gamma  = \partial\Gamma^{n-2,n-1} +\bar\partial\Gamma^{n-1,n-2} \hspace{3ex} \mbox{on}\hspace{1ex} \widetilde{X},\end{eqnarray} where $\Gamma$ is an $\widetilde\omega$-bounded $C^\infty$ $(2n-3)$-form on $\widetilde{X}$, $\Gamma^{n-2,n-1}$ (resp.$\Gamma^{n-1,n-2}$) is the $(n-2, n-1)$-component (resp. $(n-1, n-2)$-component) of $\Gamma$  and $\widetilde\omega=\pi_X^\star\omega$ is the lift of the metric $\omega$ to $\widetilde{X}$.
It is a well-established fact (see [AI01] or [IP12] for a proof or [Pop15] for a short proof) that a metric
which is both SKT and balanced is K\"ahler. Using (\ref{skt_hyper_1}), we get \begin{eqnarray}\label{SKT_balanced_1}
    \pi_X^\star\omega^{n-1} =\widetilde{\omega}^{n-1}=\widetilde{\omega}^{n-2}\wedge(\partial\alpha+\bar\partial\beta)=\partial(\widetilde{\omega}^{n-2}\wedge\alpha)+\bar\partial(\widetilde{\omega}^{n-2}\wedge\beta)
\end{eqnarray}
where the last equality is obtined by the fact that $\omega$ is a K\"ahler metric so $\widetilde{\omega}$.
Putting together (\ref{balanced_hyper_1}) and (\ref{SKT_balanced_1}), we get: $ \Gamma^{n-2,n-1} = (\widetilde{\omega}^{n-2}\wedge\alpha)$  such that:  \begin{eqnarray}\label{7_in_thm}\bar\partial\Gamma^{n-2,n-1}=  \bar\partial (\widetilde{\omega}^{n-2}\wedge\alpha)=\widetilde{\omega}^{n-2}\wedge\bar\partial\alpha \,\in \mbox{L}^2(\widetilde X)
\end{eqnarray}
and $ \Gamma^{n-1,n-2} = (\widetilde{\omega}^{n-2}\wedge\beta) $ such that: \begin{eqnarray}\label{8_in_thm}\partial\Gamma^{n-1,n-2}=\partial (\widetilde{\omega}^{n-2}\wedge\beta)=\widetilde{\omega}^{n-2}\wedge\partial\beta\,\in \mbox{L}^2(\widetilde X).
\end{eqnarray}
Since $\widetilde\omega^{n-2}$ is bounded on $\widetilde X$ with respect to itself, we get by (\ref{7_in_thm}) and (\ref{8_in_thm}) $\bar\partial\alpha\in \mbox{L}^2(\widetilde X)$ and $\partial\beta\in \mbox{L}^2(\widetilde X)$. Proposition \ref{Equi_SKT_Kahler} gives us what we want for all $n>2$. The case of a compact complex surface is trivial, as it combines both balanced hyperbolicity and Kähler hyperbolicity.
\hfill $\Box$

\vspace{2ex}

Now we will give a new proof of the implication \noindent {\it SKT hyperbolicity $\implies$ Kobayashi hyperbolicity} as an application of the following theorem
\begin{The}([Duv08]).
    A compact complex manifold is Kobayashi hyperbolic if and only if any "holomorphic disc" in $X$ satisfy a linear isoperimetric inequality.

\end{The}
This means that there exists a real constant $C > 0$ such that given any holomorphic map $f:\mathbb{D}\to X$ which is say continuously differentiable up to
the boundary, we have

$$\mbox{Area} (\mathbb{D})\leq C \mbox{Length}(\partial\mathbb{D}),$$

with respect to some (hence any, by compactness, by properly modifying
the constant $C$) Hermitian metric on $X$.
Using this theorem, it is almost immediate to prove Kobayashi hyperbolicity of SKT hyperbolic manifolds.
\begin{Cor}
    Let $X$ be an SKT hyperbolic manifold. Then, $X$ is Kobayashi hyperbolic.
\end{Cor}
\noindent {\it Proof.}
Let $X$ be a compact complex manifold, with $\mbox{dim}_\C X=n$, equipped with a {\it SKT hyperbolic} metric $\omega$. This means that, if $\pi_X:\widetilde{X}\longrightarrow X$ is the universal cover of $X$, we have $$\pi_X^\star\omega = \partial\alpha+\bar\partial\beta  \hspace{3ex} \mbox{on}\hspace{1ex} \widetilde{X},$$ where $\alpha$ and $\beta$ are a smooth, $\widetilde\omega$-bounded $(0,1)$-form and $(1,0)$-form  on $\widetilde{X}$ and $\widetilde\omega=\pi_X^\star\omega$ is the lift of the metric $\omega$ to $\widetilde{X}$. Taking any holomorphic lifting $\widetilde{f}$ of $f$ to $\widetilde{X}$, namely a holomorphic map $\tilde{f}:\mathbb{D}\longrightarrow\widetilde{X}$ such that $f=\pi_X\circ\tilde{f}$, we get the following 
 \begin{Claim}\label{Claim:f-tilde-star_gamma_bounded} The $1$-form $\tilde{f}^\star(\alpha+\beta)$ is $(f^\star\omega)$-bounded on $T_{\overline{\mathbb{D}}}$.

  \end{Claim}

  \noindent {\it Proof of Claim.} For any tangent vector $v\in T_{\overline{\mathbb{D}}}$, we have: \begin{eqnarray*}\big|\big(\tilde{f}^\star(\alpha+\beta)\big)(v)\big|^2 & = & \big|\big(\alpha+\beta\big)(\tilde{f}_\star v)\big|^2 \stackrel{(a)}{\leq} (C_1+C_2)\,|\tilde{f}_\star v|^2_{\widetilde\omega} \\
    & = & C\,|v|^2_{\tilde{f}^\star\widetilde\omega} \stackrel{(b)}{=} C\,|v|^2_{f^\star\omega},\end{eqnarray*} where $C>0$ is a constant independent of  $v$ that exists such that inequality (a) holds thanks to the {\it $\widetilde\omega$-boundedness} of $\alpha$ and $\beta$ on $\widetilde{X}$, while (b) follows from $\tilde{f}^\star\widetilde\omega = f^\star\omega$.   \hfill $\Box$
    
Now observe that
\begin{eqnarray*}
    \mbox{Area}(\mathbb{D})&=&\int_\mathbb{D}f^\star(\omega)=\int_\mathbb{D}(\pi\circ\widetilde {f})^\star\omega\\
    &=& \int_\mathbb{D}(\widetilde {f})^\star(\partial\alpha+\bar\partial\beta) \stackrel{(i)}{=}
\int_\mathbb{D}\widetilde {f}^\star d(\alpha+\beta)=\int_{\partial\mathbb{D}}\widetilde {f}^\star(\alpha+\beta)\\
&\stackrel{(ii)}{\leq}& C \mbox{Length}(\partial\mathbb{D})
\end{eqnarray*}

 where $C>0$ is a constant that exists such that inequality (ii) holds thanks to Claim \ref{Claim:f-tilde-star_gamma_bounded} and equality (i) holds since the (1,1)-form $f^\star\omega$ already has maximal degree on $\mathbb{D}$.
\hfill $\Box$ 

\begin{Prop}
    Let $(X,\omega)$, be a compact complex manifold with $\mbox{dim}_\C X=n$, and consider the  universal cover $\pi_X:\widetilde{X}\longrightarrow X$  of $X$. If  $\omega$ is both SKT hyperbolic and balanced hyperbolic
Then, there exists a constant $C > 0$ such that for any bounded domain with $C^1$ boundary $\Omega\subset \widetilde X$ we have the following linear isoperimetric inequality

$$\mbox{Volume}_{\pi^\star(\omega)}(\Omega) \leq C \mbox{Area}_{\pi^\star(\omega)}(\partial\Omega).$$
\end{Prop}
\noindent {\it Proof.}
To obtain the linear isoperimetric inequality, it suffices to consider the case where $\omega$ is simultaneously balanced hyperbolic and  SKT (similarly, if $\omega$ is assumed to be SKT hyperbolic and balanced). Let consider the case when $\omega$ is both balanced hyperbolic and SKT, so $\omega$ is also a K\"ahler metric, we get 
\begin{eqnarray*}
\mbox{Volume}_{\pi^\star(\omega)}(\Omega)=\int_\Omega\pi^\star\omega^n = \int_\Omega\pi^\star\omega^{n-1}\wedge\pi^\star\omega&=&\int_\Omega d(\Gamma\wedge\pi^\star\omega)\\&=&\int_{\partial\Omega} \Gamma\wedge\pi^\star\omega\\ &\leq& C \mbox{Area}_{\pi^\star(\omega)}(\partial\Omega),
\end{eqnarray*}
where the last equality is derived from Stokes's theorem, while the final inequality is obtained due to the boundedness of $\Gamma$.
\hfill $\Box$

\vspace{3ex}
\noindent {\bf References.} \\
\noindent [Bro78]\, R. Brody --- {\it Compact Manifolds and Hyperbolicity} --- Trans. Amer. Math. Soc. {\bf 235} (1978), 213-219.  

\vspace{1ex}

\noindent [Dem97]\, J.-P. Demailly --- {\it Complex Analytic and Algebraic Geometry} --- http://www-fourier.ujf-grenoble.fr/~demailly/books.html
\vspace{1ex}


\noindent [DE00] Demailly, J.-P., El Goul, --- {\it  Hyperbolicity of Generic
Surfaces of High Degree in Projective 3-Space} --- Amer. J. Math. 122 (2000), no. 3, 515–546.
\vspace{1ex}

\noindent [Gro91]\, M. Gromov --- {\it K\"ahler Hyperbolicity and $L^2$ Hodge Theory} --- J. Diff. Geom. {\bf 33} (1991), 263-292.

\vspace{1ex}
\noindent [HL83]\, R. Harvey, H.B. Lawson --- {\it An intrinsic characterization of K¨ahler manifolds} --- Invent. Math. {\bf 74} (1983), no. 2, 169–198.

\vspace{1ex}

\noindent [Kob70]\, S. Kobayashi --- {\it Hyperbolic Manifolds and Holomorphic Mappings} --- Marcel Dekker, New York (1970).








\vspace{1ex}
\noindent [MP22a]\, S. Marouani, D. Popovici --- {\it Balanced Hyperbolic and Divisorially Hyperbolic Compact Complex Manifolds } ---arXiv e-print CV 2107.08972v2, to appear in Mathematical Research Letters.

\vspace{1ex}

\noindent [MP22b]\, S. Marouani, D. Popovici --- {\it Some properties of balanced hyperbolic
compact complex manifolds } --- Internat. J. Math. 33 (2022), no. 3, Paper No.
2250019, 39. MR 4390652.


\vspace{1ex}

\noindent [Pop15]\, D. Popovici --- {\it Aeppli Cohomology Classes Associated with Gauduchon Metrics on Compact Complex Manifolds} --- Bull. Soc. Math. France {\bf 143}, no. 4 (2015), p. 763-800.



\vspace{1ex}

\noindent [Sch07]\, M. Schweitzer --- {\it Autour de la cohomologie de Bott-Chern} --- arXiv e-print math.AG/0709.3528v1.
\vspace{1ex}
\noindent [YZZ22]\,S.T. Yau, Q. Zhao, F. Zheng, On Strominger K\"ahler-like manifolds with degenerate torsion,
preprint arXiv:1908.05322 [math.DG].
\vspace{1ex}

\noindent [ST10]\, J. Streets, G. Tian --- {\it A Parabolic Flow of Pluriclosed Metrics} --- Int. Math. Res. Notices, {\bf 16} (2010), 3101-3133.

\vspace{1ex}

\noindent [Voi02]\, C. Voisin --- {\it Hodge Theory and Complex Algebraic Geometry. I.} --- Cambridge Studies in Advanced Mathematics, 76, Cambridge University Press, Cambridge, 2002.



\vspace{1ex}

\noindent [YZZ23]\, S.-T. Yau, Q. Zhao, F. Zheng --- {\it On Strominger K\"ahler-like Manifolds with Degenerate Torsion} --- arXiv e-print DG 1908.05322v2

\vspace{3ex}

\noindent Institut de Math\'ematiques de Toulouse, Universit\'e Paul Sabatier,

\noindent 118 route de Narbonne, 31062 Toulouse, France

\noindent Email: 	almarouanisamir@gmail.com

\end{document}